\documentclass{amsproc}
\usepackage{amssymb,amsmath,amscd}

\def\N{{\mathbb N}}
\def\Z{{\mathbb Z}}

\def\R{{\mathbb R}}

\def\lin{{\rm lin}}
\def\Vol{{{\rm Vol}}} 
\def\vol{{\rm vol}}
\def\conv{{\rm conv}}

\def\deg{{\rm deg}}
\def\st{{\rm st}}

\newtheorem{theorem}{Theorem}[section]
\newtheorem{lemma}[theorem]{Lemma}

\theoremstyle{definition}

\newtheorem{coro}[theorem]{Corollary}
\newtheorem{conj}[theorem]{Conjecture}

\theoremstyle{remark}
\newtheorem{rem}[theorem]{Remark}

\numberwithin{equation}{section}

\begin{document}

\title{Lattice Polytopes of Degree $2$}

\author{Jaron Treutlein}
\address{Department of Mathematics and Physics, 
University of T\"ubingen, Auf der Morgenstelle 10, D-72076 T\"ubingen, Germany}
\email{jaron@mail.mathematik.uni-tuebingen.de}





\begin{abstract}
A theorem of Scott gives an upper bound for the normalized volume of 
lattice polygons with exactly $i>0$ interior lattice points. We will 
show that the same bound is true for the normalized volume of lattice 
polytopes of degree 2 even in higher dimensions. In particular, there 
is only a finite number of quadratic polynomials with fixed leading 
coefficient being the $h^*$-polynomial of a lattice polytope.
\end{abstract}

\maketitle

\section{Introduction}

An $n$-dimensional lattice polytope $P\subset \R^n$ is the convex hull 
of a finite number of elements of $\Z^n$. In the following, we denote by 
$\Vol(P)=n!\vol(P)$ the normalized volume of $P$ and may call it the 
volume of $P$. By $\Pi^{(1)}:=\Pi(P)\subset\R^{n+1}$, we denote the 
convex hull of $(P,0)\subset\R^{n+1}$ and $(0,\ldots,0,1)\in\R^{n+1}$, 
which we will call the standard pyramid over $P$. Recursively we define 
$\Pi^{(k)}(P)=\Pi\Big(\Pi^{(k-1)}(P)\Big)$ for all $k>0$. $\Delta_n$
will denote the $n$-dimensional basic lattice simplex throughout, i.e. 
$\Vol(\Delta_n)=1$. If two lattice polytopes $P$ and $Q$ of the same 
dimension are equivalent via some affine unimodular transformation, we will 
write $P\cong Q$. The $k$-fold of a polytope $P$ will be the convex hull of
the $k$-fold vertices of $P$ for every $k\geq0$.

\medskip

Pick's formula gives a relation between the normalized volume, the 
number of interior lattice points and the number of lattice points of 
a lattice polygon, i.e. of a two-dimensional lattice polytope: $\Vol(P)
=|P\cap\Z^2|+|P^\circ\cap\Z^2|-2$. Here $P^\circ$ means the interior of
the polytope $P$.

\medskip

In 1976 Paul Scott \cite{Sc} proved that the volume of a lattice 
polygon with exactly $i\geq1$ interior lattice points is constrained 
by $i$:

\begin{theorem}[Scott]
Let $P\subset\R^2$ be a lattice polygon such that 
$|P^\circ\cap\Z^2|$ $=i\geq1$. If $P\cong 3\Delta_2$, then $\Vol(P)=9$ 
and $i=1$. Otherwise the normalized volume is bounded 
by $\Vol(P)\leq 4(i+1)$. According to Pick's formula, this implies 
$|P\cap\Z^2|\leq 3i+6$ and $|P\cap\Z^2|\leq\frac{3}{4}\Vol(P)+3$.
\label{Scott}
\end{theorem}

Besides Scott's proof, there are two proofs by Christian Haase and 
Joseph Schicho \cite{HS}. Another proof is given in \cite{Tr}.

\medskip

Our aim is to generalize Scott's theorem. Therefore we need to 
introduce another invariant, the degree of a lattice polytope:\\
It is known from \cite{Eh}, \cite{St1} and \cite{St2} that $h^*_P(t) := 
(1-t)^{n+1}\sum_{k\geq0}|kP \cap \Z^n|t^k\in\Z[t]$ is a polynomial of 
degree $d\leq n$. This number is described as the degree of $P$ and is 
the largest number $k\in\N$ such that there is an interior lattice 
point in $(n+1-k)P$ (cf. \cite{BN}). The leading coefficient of 
$h^*_P$ is the number of interior lattice points in $(n+1-d)P$ and the 
constant coefficient is $h_P^*(0)=1$. Moreover the sum of all coefficients
is the normalized volume of $P$ and all coefficients are non-negative 
integers by the non-negativity theorem of Richard P. Stanley \cite{St1}.\\
It is easy to show that the $h^*$-polynomial of $P$ and $\Pi(P)$ are 
equal. So $P$ and $\Pi(P)$ have the same degree and the same normalized 
volume, which is the sum of all coefficients of the $h^*$-polynomial. 
Moreover \[\Big|\Big((n+2-d)\Pi(P)\Big)^\circ\cap\Z^{n+1}\Big|=\Big|
\Big((n+1-d)P\Big)^\circ\cap\Z^n\Big|.\]
Scott's theorem shows that the normalized volume of a two-dimensional 
lattice polytope of degree 2 with exactly $i>0$ interior 
lattice points is bounded by $4(i+1)$, except for one single polytope: 
$3\Delta_2$. We generalize this result to the case of $n$-dimensional 
lattice polytopes of degree 2.

\begin{theorem}
Let $P\subset\R^n$ be a $n$-dimensional lattice polytope of degree $2$.
If $P\cong\Pi^{(n-2)}(3\Delta_2)$, then $\Vol(P)=9,$ $|P\cap\Z^n|=8+n$ 
and $\Big|\Big((n-1)P\Big)^\circ\cap\Z^n\Big|=1$. Otherwise the 
following equivalent statements hold: \begin{align*} (1)&\ \Vol(P)\leq 
4(i+1)\\
(2)&\ b\leq 3i+n+4\\(3)&\ b\leq\frac{3}{4}\Vol(P)+n+1,\end{align*} 
where $b:=|P\cap\Z^n|$ and  $i:=\Big|\Big((n-1)P\Big)^\circ\cap\Z^n
\Big|\geq1$.
\label{Deg2}
\end{theorem}

The following theorem of Victor Batyrev \cite{Ba} motivates our estimation 
of the normalized volume of a lattice polytope of degree $d$:

\begin{theorem}[Batyrev] 
Let $P\subset\R^n$ be an $n$-dimensional lattice polytope of degree $d$. 
If \[n\geq 4d \binom{2d+\Vol(P)-1} {2d},\] then $P$ is a standard pyramid 
over an $(n-1)$-dimensional lattice polytope.
\label{finit}
\end{theorem}

There is a recent result by Benjamin Nill \cite{Ni} which even 
strenghtens this bound:

\begin{theorem}[Nill]
Let $P \subset \R^n$ be a $n$-dimensional lattice polytope of degree 
$d$. If \[n\geq (\Vol(P)-1)(2d+1),\]
then $P$ is a standard pyramid over an $(n-1)$-dimensional lattice 
polytope.
\label{Benjamin}
\end{theorem}

Jeffrey C. Lagarias and G\"unter M. Ziegler showed in \cite{LZ} that up to 
unimodular transformation there is only a finite number of 
$n$-dimensional lattice polytopes having a fixed volume. From Theorem 
\ref{finit} or Theorem \ref{Benjamin} follows
\begin{coro}[Batyrev]
For a family $\mathcal{F}$ of lattice polytopes of degree $d$, the 
following is equivalent:
\begin{enumerate}\item[$(1)$] $\mathcal{F}$ is finite modulo standard pyramids 
and affine unimodular transformation,
\item[$(2)$] There is a constant $C_d>0$ such that $\Vol(P)\leq C_d$ for all 
$P\in \mathcal{F}$.
\end{enumerate}
\end{coro}
\begin{conj}[Batyrev]
Let $P$ be a lattice polytope of degree $d$ with exactly $i\geq1$ 
interior lattice points in its $(\dim(P)+1-d)$-fold. Its normalized 
volume $\Vol(P)$ can then be bounded by a constant $C_{d,i},$ only 
depending on $d$ and $i$. The finiteness of lattice polytopes of degree 
$d$ with this property up to standard pyramids and affine unimodular 
transformation follows from Theorem \ref{finit}.
\label{Verm}
\end{conj}
 Theorem \ref{Deg2} proves Conjecture \ref{Verm} in the case $d=2$.

\noindent
\begin{coro}
Up to affine unimodular transfomations and standard pyramids there is 
only a finite number of lattice polytopes of degree $2$ having exactly 
$i\geq 1$ interior lattice points in their adequate multiple.
\label{Oma}
\end{coro}

\noindent
This follows from Theorem \ref{Deg2} and Theorem \ref{finit}.
\begin{coro} There is only a finite number of quadratic polynomials 
$h\in\Z[t]$ with leading coefficient $i\in\N$, such that $h$ is the 
$h^*$-polynomial of a lattice polytope.
\end{coro}

\noindent
This follows from Theorem \ref{Deg2} and the fact 
that all coefficients of $h^*_P$ are positive integers summming up to 
$\Vol(P)$.\\

\medskip

In the remaining part of the paper we prove Theorem \ref{Deg2}.\\

\medskip
\noindent {\bf Acknowledgments:} The author would like to thank Victor 
Batyrev and Benjamin Nill for discussions and joint work on this subject.

\section{Preparations}

The formula of Pick can be easily generalized for higher dimensional 
polytopes of degree $2$ using their $h^*$-polynomial. This shows that
statements (1) -- (3) in Theorem \ref{Deg2} are equivalent.
\begin{lemma} 
An $n$-dimensional lattice polytope of degree $2$ has normalized volume 
$\Vol(P)=b+i-n,$ where $b:=|P\cap\Z^n|$ and $i:=\Big|\Big((n-1)P
\Big)^\circ\cap\Z^n\Big|$.
\label{VPick} 
\end{lemma}

\noindent
\begin{proof}{Proof.} The normalized volume of $P$ can be computed by 
adding the coef- ficients of the $h^*$-polynomial of $P$. Consequently 
$\Vol(P)= 1+(b-n-1)+i.$
\end{proof}

Let $s\subset P$ be a face of $P$. By $\st(s)=\bigcup F,$ we denote the 
star of $s$ in $P$, where the union is over all faces $F\subset P$ 
of $P$ containing $s$.
\begin{lemma}
Let $P$ be an $n$-dimensional lattice polytope of degree $2$ and 
$s\subset P$ a face of $P$ having exactly $j>0$ interior lattice 
points in its $(n-2)$-fold: \[\Big((n-2)s\Big)^\circ\cap\Z^n=
\{x_1,\ldots,x_j\}.\]
Moreover, we suppose \[z:=\Big|P\backslash\st(s)\cap\Z^n\Big|\geq 1.\]
Then $0<j+z-1\leq \Big|\Big((n-1)P\Big)^\circ\cap\Z^n\Big|$.
\label{lines}
\end{lemma}

\begin{rem}
{\rm Let us first consider an easy case.\\
 If $z=1$, i.e. $P\backslash\st(s)\cap\Z^n=\{p\}$, then  \[p+x=(n-1)\Big(\frac{n-2}{n-1}\frac{x}{n-2}+\frac{p}{n-1}\Big)\in\Big(
(n-1)P\Big)^\circ\cap\Z^n \ \forall x\in\Big((n-2)s\Big)^\circ\cap\Z^n\]
yield $j>0$ distinct lattice points in $(n-1)P$. So $0<j\leq \Big|
\Big((n-1)P\Big)^\circ\cap\Z^n\Big|$ as claimed.}
\label{z=1}
\end{rem}

\begin{proof}
If $l=1$, the claim is certainly correct. Hence let $l\geq2$.\\
There is a lattice point $z_l\in y_l^\perp\cap\Big(D\backslash\{
\frac{x_1}{n-2}\}\Big)\cap\Z^{n+1}$. Define $\pi_l:=\conv(s_l,z_l)$. 
Obviously $\pi_l\cap s=s_l$. By induction, there are further pyramids 
$\pi_1,\ldots,\pi_{l-1}$ satisfying $\pi_k\cap s=s_k$ and $\pi_k\cap\pi_{k'}\subset\{z_1,\ldots,z_{k'}\}\subset\partial 
D\cap\Z^{n+1} \ \forall k<k'<l$.

\smallskip

Assume $\pi_l\cap\pi_k \not\subset\{z_1,\ldots,z_l\}$, i.e. there 
exists a point $q\in\pi_l\cap\pi_k$, $q\not\in\{z_1,\ldots,z_l\}$ and 
$k<l$. Therefore $y_k|_{\pi_k}\geq 0$, because $y_k|_{s_k}\geq 0$ and 
$y_k(z_k)=0$. In particular, $y_k(q)\geq0$, and $y_l(q)\geq0$ as well. 
As $q\in\pi_l=\conv(s_l,z_l)$, there is a point $p\in s_l$ and a number 
$\lambda\in[0,1]$ such that $q=\lambda p+(1-\lambda)z_l$. Therefore 
$0\leq y_k(q)=\lambda y_k(p)+(1-\lambda)y_k(z_l)$ with $y_k(z_l)\leq0$ 
as $z_l\in D$ and $y_k|_D\leq 0$.

\smallskip

If $y_k(p)\geq0$, then $p\in s\cap\{x\in\R^{n+1}\ : \ y_k(x)\geq0\}=
s_k'\subseteq\bigcup_{r\leq k}s_r$. But this is a contradiction to 
$p\in s_l$ with $l>k$. So \[0\leq y_k(q)=\lambda y_k(p)+(1-\lambda)
y_k(z_l)\leq0\] with equality only in the case of $\lambda=0$ and 
$y_k(z_l)=0$. Therefore the intersection of $\pi_k$ and $\pi_l$ is 
$q=z_l$ or empty. This is a contradiction to $\pi_l\cap\pi_k \not
\subset\{z_1,\ldots,z_l\}$ and so the claim is proven.\hfill $\Box$

\medskip

%
%
The pyramids $\pi_1,\ldots,\pi_k$ intersect with $D$ only in faces of 
$D$.\\
To any $k\in\{1,\ldots, K\}$ denote by $a_k:=\Big|\Big((n-2)s_k
\Big)^\circ\cap\Z^{n+1}\Big|$ the number of interior lattice points of 
$(n-2)s_k$. By Remark \ref{z=1}, there are $a_k\geq0$ interior lattice 
points of $(n-1)s$ in $(n-1)\pi_k$. By adding up the number of 
interior lattice points in $(n-1)\pi_1,\ldots,(n-1)\pi_K$, we derive 
from the claim 
\[\Big|\bigcup_{k=1}^K\Big((n-1)\pi_k\Big)^\circ\cap\Z^{n+1}\Big|
\geq\sum_{k=1}^K a_k=j-1.\]
Furthermore to every $p\in D\backslash \{\frac{x_1}{n-2}\}$ we get a 
lattice point of $\Big((n-1)\Big(D\backslash\{\frac{x_1}{n-2}\}\Big)
\Big)^\circ$ $\subset \Big((n-1)P\Big)^\circ$ in the following way: \[p+x_1=(n-1)\Big(\frac{n-2}{n-1}\frac{x_1}{n-2}+\frac{p}{n-1}\Big)
\in\Big((n-1)D\Big)^\circ\cap\Z^{n+1}.\]
Finally we get $\Big|\Big((n-1)P\Big)^\circ\cap\Z^{n+1}\Big|\geq 
j-1+z$.\hfill
\end{proof}

\section{The Proof of the Main Theorem}

\noindent
If $n=2$, then Theorem \ref{Deg2} is equal to Scott's Theorem 
\ref{Scott}. So let $n>2$.\\
The monotonicity theorem of Stanley \cite{St3} says that the degree 
of every face of a polytope is not greater than the degree of the 
polytope itself. In particual this is true for every facet. So we will
distinguish the two cases that there is a facet of $P$ having degree $2$
or there is not.\\ 
For the second case we need a result of Victor Batyrev and Benjamin 
Nill. They proved in \cite{BN} that every $n$-dimensional lattice 
polytope of degree less than 2 either is equivalent to a pyramid over 
the exceptional lattice simplex $2\Delta_2$ or it is a Lawrence polytope, 
i.e. a lattice polytope projecting along an edge onto an 
$(n-1)$-dimensional basic simplex.

\medskip

 \textbf{Case 1:} There is a facet $F\subset P$ of $P$ 
having degree two, i.e. \[\Big|\Big((n-2)F\Big)^\circ\cap\Z^n\Big|
=j\geq1.\]
Define $z:=|P\backslash F\cap\Z^n|$. From Lemma \ref{lines} we get 
$z+j-1\leq i$. Thus, by induction, we get, if $F\not\cong
\Pi^{(n-3)}(3\Delta_2)$,
\begin{align*}|P\cap\Z^n|&=|F\cap\Z^n|+|(P\backslash F)\cap\Z^n|\leq 
3j+n-1+4+z\\
&=3(j+z-1)-2z+2+n+4\stackrel{z\geq1}{\leq}3i+n+4,\end{align*}
\noindent Otherwise $F\cong\Pi^{(n-3)}(3\Delta_2)$ and again by 
induction and Lemma \ref{lines}:\\ $|F\cap\Z^n|=(n-1)+8$, $z\leq i$ 
and so  $|P\cap\Z^n|=n-1+8+z\leq i+7+n$. This term is smaller than 
$3i+n+4$ if $i\geq2$. If $i=1$ however, we get \[n+8\leq|P\cap\Z^n|
=n+7+z\leq i+7+n=8+n,\] so $|P\cap\Z^n|=8+n$ and $\Vol(P)=9$ by Lemma
\ref{VPick}. In this case $P\cong\Pi^{(n-2)}(3\Delta_2)$ because 
$\Vol(F)=9$ and $F\cong\Pi^{(n-3)}(3\Delta_2)$.\\
%
%
%

\medskip

\textbf{Case 2:} Every facet $F$ of $P$ has degree $\deg(F)\leq 1.$

\smallskip

Let $y$ be an edge of $P$ having the maximal number of lattice points; 
its length will be denoted by $h_1$, i.e. $h_1=|y\cap\Z^n|-1$. Among 
all $2$-codimensional faces of $P$ containing $y$, $s$ 
should be the face having the maximal number of lattice points. We will
denote by $F_1$ and $F_2$ the two facets of $P$ containing $s$.

\smallskip 

 Again the monotonicity theorem of Stanley \cite{St3} implies $\deg(s)
\leq \deg(F_1)=1$. Similarly to case 1, we will denote by $z:=|P
\backslash\{F_1\cup F_2\}\cap\Z^n|$ the number of lattice points of 
$P$ not in $F_1$ and $F_2$.

\smallskip

By the result of Victor Batyrev and Benjamin Nill \cite{BN} we find
that the facets $F_1$ and $F_2$ are either $(n-1)$-dimensional Lawrence 
polytopes or pyramids over $2\Delta_2$.\\
%
%
%

\medskip

\textbf{(A)} $F_1$ and $F_2$ are Lawrence polytope with heights $h_1^{(k)}, 
h_2^{(k)},\ldots, h_{n-1}^{(k)}\ \forall  k\in\{1,2\}$, where we assume that
$h_l^{(1)}=h_l^{(2)}=h_l \ \forall l\in\{1,\ldots,n-2\}$,
\[s=\conv(0,h_1e_1,e_l,e_l+h_le_1 \ : \ 2\leq l
\leq n-2),\] where $\{e_1,\ldots,e_{n-2},e_{n-1}^{(k)}\}$ should denote a
lattice basis of $\lin(F_k)\cap \Z^n$ such that $F_k=\conv(s,e_{n-1}^{(k)},
e_{n-1}^{(k)}+h_{n-1}^{(k)}e_1)$ for $k\in\{1,2\}$. Since the degree of the 
Lawrence prism $s$ is at most one, we obtain \[\Big|\Big((n-2)s\Big)^\circ
\cap\Z^n\Big|=\Vol(s)-1=\Big(\sum_{l=1}^{n-2}h_l\Big) -1.\]
We may assume $z=|(P\backslash\{F_1\cup F_2\})\cap\Z^n|\not=0$ because 
otherwise $P$ would be a prism over the face $P\cap\{X_1=0\}$, which is 
an $(n-1)$-dimensional lattice simplex of degree at most $1$, whose only 
lattice points are vertices. By \cite{BN} this is a basic simplex and hence
$P$ is a Lawrence polytope. Consequently $\deg(P)<2$, a contradiction.\\
%
%
%
We have to distinguish the following two cases: 

\medskip

\textbf{(i)} $\Big|\Big((n-2)s\Big)^\circ\cap\Z^n\Big|\geq1$. 

\smallskip

Because of Lemma \ref{lines}, we get the estimation \[z+\Big(\Big(
\sum_{l=1}^{n-2}h_l\Big)-1\Big)-1\leq i.\]
So we can bound the number of lattice points of $P$:
\begin{eqnarray*}|P\cap\Z^n|&=&|(F_1\cup F_2)\cap\Z^n|+z =|s\cap\Z^n|+
h_{n-1}^{(1)}+1+h_{n-1}^{(2)}+1+z\\
&=&\sum_{l=1}^{n-2} h_l + (n-2)+ h_{n-1}^{(1)}+h_{n-1}^{(2)}+2+z\leq 
i+n+2h_1+2\\
&\stackrel{h_1\leq i+1}{\leq}&i+n+2(i+1)+2=3i+n+4.\end{eqnarray*}
%
%

\medskip

\textbf{(ii)} $\Big|\Big((n-2)s\Big)^\circ\cap\Z^n\Big|=0$. 

\smallskip

In this case, $s$ has degree zero, so it is a basic simplex. Our 
assumption on $s$ implies that every lattice point of $P$ is a vertex.
If $n=3$, then Howe's theorem \cite{Sca} yields that $P$ has at most 
$8$ vertices, therefore $|P\cap\Z^n|\leq8<n+4+3i$. So let $n\geq4$. 

\smallskip

In that case, since every $2$-codimensional face is a simplex and every facet 
is a Lawrence prism, we see that $P$ is simplicial, i.e. every facet 
is a simplex. We may suppose that $P$ is not a simplex. Let $S$ be a 
subset of the vertices of $P$ such that the convex hull of $S$ is not 
a face of $P$. Then the sum over the vertices of $S$ is a lattice point 
in the interior of $|S|\cdot P$. Since the degree of $P$ is two, this 
implies $|S|\geq n-1$. In other words, every subset of the vertices of 
$P$ that has cardinality at most $n-2$ forms the vertex set of a face of 
$P$, i.e. $P$ is $(n-2)$-neighbourly. As is known from 
\cite{Bro}, a polytope of dimension $n$ that is not a simplex is at 
most $\lfloor \frac{n}{2}\rfloor$-neighbourly. Therefore $n-2\leq
\frac{n}{2}$. This shows $n=4$. 

\smallskip

Let $f_j\geq0$ be the number of $j$-dimensional faces of $P$. 
Since $P$ is a $2$-neighbourly simplicial $4$-dimensional polytope 
we get $f_1={f_0\choose2}$ and $f_2=2f_3$. Since the Euler 
characteristic of the boundary of $P$ vanishes, i.e. $f_0-f_1+f_2-f_3=0$, 
we deduce $f_3=\frac{f_0(f_0-3)}{2}$. Let $\mathcal D$ denote the set of 
subsets $\Delta$ of the vertices of $P$ such that $\Delta$ has cardinality 
three but $\Delta$ is not the vertex set of a face of $P$. Therefore, 
$|\mathcal D|$ $={f_0\choose 3}-f_2=f_0\Big(\frac{(f_0-1)(f_0-2)}{6}-(f_0-3)\Big)$. 
Since $|\{(e,\Delta)\ : \ e \mbox{ is an edge of }P,\Delta\in\mathcal D,$
$e\subset\mathcal D\}|$ $=3|\mathcal D|$, double counting yields that there exists 
an edge $e$ of $P$ that is contained in at least $\frac{3|\mathcal D|}{f_1}$ many 
elements $\Delta\in\mathcal D$. Therefore, any such $\Delta$ contains one vertex 
that is not in the star of $e$, and hence Lemma \ref{lines} yields
\[i\geq\frac{3|{\mathcal D}|}{f_1}=f_0-2-6\frac{f_0-3}{f_0-1}\geq f_0-8.\]
Thus, $|P\cap\Z^n|=f_0\leq 8+i<n+4+3i$.

\medskip

%
%
%
\textbf{(A')} $F_1$, $F_2$ and $s$ have no common projection direction.\\
Without loss of generality let $F_1$ and $s$ have two different projection
directions. If $s$ contains an edge of length at least $2$, then this has 
to be a common projection direction with $F_1$, because $s$ and $F_1$ are 
Lawrence prisms. But this is a contradiction. Hence, all lattice points in 
$s$ are vertices. In particular, $y$ has length one, so also all lattice 
points of $P$ are vertices. \\
Since any of the two different projection directions of the Lawrence prism 
$s$ maps a four-gon face onto the edge of an unimodular base simplex 
and two edges of the four-gon give the projection direction, we see that 
there is at most one four-gon face in $s$. Therefore, $s$ contains at most 
$(n-2)+2=n$ lattice points.\\
Since $F_k$ contains at most two vertices not in $s$ for $k\in\{1,2\}$, we 
get $|(F_1\cup F_2)\cap\Z^n|\leq n+4<n+4+3i$. Therefore we may assume 
$z:=|P\backslash(F_1\cup F_2)\cap\Z^n|\not=0$.\\
If $\Big|\Big((n-2)s\Big)^\circ\cap\Z^n\Big|=0$, then we will proceed exactly
like in case (ii) from (A). So let 
$j:=\Big|\Big((n-2)s\Big)^\circ\cap\Z^n\Big|\geq1$. \\
Because of Lemma \ref{lines}, we get the estimation $z+j-1\leq i,$ in 
particular $z\leq i$. Hence we can bound the number of lattice points 
of $P$:
\begin{eqnarray*}|P\cap\Z^n|=|(F_1\cup F_2)\cap\Z^n|+z \leq n+4+i<3i+n+4.
\end{eqnarray*}
%
%
%

\medskip

\textbf{(B)} $F_1$ is a Lawrence polytope with the heights $h_1\geq 
h_2\geq\ldots\geq h_{n-1}$, $F_2\cong\Pi^{(n-3)}(2\Delta_2)$. \\
Here \[s\cong\conv(0,h_1e_1,e_l, \ 2\leq l\leq n-2)\]
and $h_1=2$, $h_2=\cdots=h_{n-2}=0$, because $s$ is contained in 
the simplex $F_2$. If $z=|P\backslash\{F_1\cup F_2\}\cap\Z^n|=0,$ 
then 
\begin{eqnarray*}
|P\cap\Z^n|&=&|F_2\cap\Z^n|+|F_1\backslash F_2\cap\Z^n|=6+(n-3)
+h_{n-1}+1\\
&\stackrel{h_{n-1}\leq h_1=2}{\leq}&4+n+2<3i+n+4.
\end{eqnarray*}
Otherwise if $z\geq1$, we obtain just like in (A) $0<z+(
h_1-1)-1\leq i$. Therefore
\begin{eqnarray*}|P\cap\Z^n|&=&|(F_1\cup F_2)\cap\Z^n|+z=|s\cap\Z^n|+
(h_{n-1}+1)+3+z\\
&=& h_1+(n-2)+(h_{n-1}+1)+3+z\leq i+4+h_{n-1}+n\\
&\stackrel{h_{n-1}\leq h_1=2}{\leq}&3i+n+4.\end{eqnarray*}
%
%
%

\medskip

\textbf{(C)} $F_1\cong F_2\cong\Pi^{(n-3)}(2\Delta_2)$.\\
Here either $s$ is a pyramid over $2\Delta_1$ or $s\cong\Pi^{(n-4)}
(2\Delta_2).$ Again $h_1=2$.\\
If $z=|P\backslash\{F_1\cup F_2\}\cap\Z^n|=0,$ 
then 
\begin{eqnarray*}
|P\cap\Z^n|=|F_2\cap\Z^n|+|F_1\backslash F_2\cap\Z^n|\leq6+(n-3)+3<3i+n+4.
\end{eqnarray*}
Otherwise if $z\geq1$, we obtain $z\leq i$ because of $\Big|\Big((n-2)s
\Big)^\circ\cap\Z^n\Big|\geq1$ and Lemma \ref{lines}. So as a result
\begin{eqnarray*}|P\cap\Z^n|&=&|F_1\cap\Z^n|+|F_2\backslash F_1\cap
\Z^n|+z\leq(6+n-3)+3+z=n+z+6\\&\leq& n+i+6\leq n+3i+4.\end{eqnarray*} 
This completes the proof.\hfill$\qed$
\pagebreak
\begin{rem} {\rm In \cite{St2}, Stanley shows that the coefficients of
$h^*_P$ also appear in the polynomial $(1-t)^{n+1}\sum_{k\geq0}|
(kP)^\circ \cap \Z^n|t^k\in\Z[t]$. So we can also compute the 
coefficients of $h^*_P$ in a different way than in Lemma \ref{VPick}.
Then it is easy to show that the bounds of 
Theorem \ref{Deg2} are also equivalent to the following estimations: \begin{eqnarray*}|(nP)^\circ\cap\Z^n|&\leq&(n+4)i+3, \\
|2P\cap\Z^n|&\leq&(4+3n)(i+1)+\frac{n(n+3)}{2}.\end{eqnarray*}}
\end{rem}
\bigskip

\bibliographystyle{amsalpha}

\end{document}